\documentclass{amsart}
\usepackage{amssymb, amsthm, amsmath}
 
\newtheorem{thm}{Theorem}[section]
\newtheorem{lem}[thm]{Lemma}
\newtheorem{cor}[thm]{Corollary}
 
\theoremstyle{definition}
\newtheorem{defn}[thm]{Definition}

\theoremstyle{remark}

\newtheorem{disc}[thm]{Discussion}
 
\numberwithin{equation}{section}

\def\CM{Cohen--Macaulay }
\def\Spec{\mathop{\mathrm{Spec}}}
\def\Cl{\mathop{\mathrm{Cl}}}

\begin{document}

\large
 
\title{$\mathbb Q$--Gorenstein splinter rings of characteristic $p$
are F--regular}

\author{Anurag K. Singh}

\address{{\em Present address:} Department of Mathematics, University of
Illinois, 1409 W. Green Street, Urbana, IL 61801, USA 
{\em e-mail:} {\tt singh6@math.uiuc.edu}}

\maketitle
 
\begin{center}

Department of Mathematics, University of Michigan, Ann Arbor, USA

{\em Received: $5$ May, $1998$,  revised: $6$ October, $1998$}
 
\end{center}

\section{Introduction}
 
A Noetherian integral domain $R$ is said to be a {\it splinter}\/ if it is a
direct summand, as an $R$--module, of every module--finite extension ring, see
\cite{Ma}. In the case that $R$ contains the field of rational numbers, it is
easily seen that $R$ is splinter if and only if it is a normal ring, but the
notion is more subtle for rings of characteristic $p>0$. It is known that
F--regular rings of characteristic $p$ are splinters, and Hochster and Huneke
showed that the converse is true for locally excellent Gorenstein rings,
\cite{HHjalg}. In this paper we extend their result by showing that
${\mathbb Q}$--Gorenstein splinters are F--regular. Our main theorem is:

\begin{thm} Let $R$ be a locally excellent ${\mathbb Q}$--Gorenstein integral
domain of characteristic $p>0$. Then $R$ is  F--regular if and only if
it is a splinter.   
\label{main}  
\end{thm}

These issues are closely related to the question of whether the {\it tight
closure} $I^*$ of an ideal $I$ of a characteristic $p$ domain agrees with its
{\it plus closure}, i.e., $I^+ = IR^+ \cap R$, where $R^+$ denotes the integral
closure of $R$ in an algebraic closure of its fraction field. We always have
the containment $I^+ \subseteq I^*$, and Smith showed that equality holds if 
$I$ is a parameter ideal in an excellent domain $R$, see \cite{Sminv}. An
excellent domain $R$ of characteristic $p$ is splinter if and only if for all
ideals $I$ of $R$, we have $I^+ =I$.

For an excellent local domain $R$ of characteristic $p$, Hochster and
Huneke showed that $R^+$ is a big \CM algebra, see \cite{HHbigcm}. For
further work on $R^+$ and plus closure see \cite{Aberbach, AH}. Our
main references for the theory of tight closure are
\cite{HHjams, HHbasec, HHjalg}.

Although tight closure is primarily a characteristic $p$ notion, it has strong
connections with the study of singularities of algebraic varieties over fields
of characteristic zero. For $\mathbb Q$--Gorenstein rings essentially of finite
type over a field of characteristic zero, it is known that F--regular type is
equivalent to log--terminal singularities, see \cite{Hara, Smfrat,
Smvanish, Walog}. Consequently our main theorem offers a
characterization of log--terminal singularities in characteristic
zero, see Corollary \ref{splinter-type}.

\section{Preliminaries}

By the {\it canonical ideal}\/ of a \CM normal domain $(R,m)$, we shall mean an
ideal of $R$ which is isomorphic to the canonical module of $R$. We next record
some results that we shall use later in our work.

\begin{lem} Let $(R,m)$ be a \CM local domain with canonical ideal $J$. Fix a
system of parameters $y_1, \dots, y_d$ for $R$, and let $s \in J$ be an element
which represents a socle generator in $J / (y_1, \dots, y_d)J$. Then for $t \in
\mathbb N$, the element $s(y_1 \dots  y_d)^{t-1}$ is a socle generator in  $J /
(y_1^t, \dots, y_d^t)J$. The ideals $I_t = (y_1^t, \dots, y_d^t)J :_R s$ form a
family of irreducible ideals which are cofinal with the powers of the maximal
ideal $m$ of $R$.  
\label{irreducible}
\end{lem}

\begin{proof} See the proof of \cite[Theorem 4.6]{HHjalg}.
\end{proof}

\begin{lem} Let $R$ be a \CM normal domain with canonical ideal $J$. Pick
$y_1 \neq 0$ in $J$. Then there exists an element $y_2$ not in any minimal prime
of $y_1$ and $\gamma \in J$ such that $y_2^i J^{(i)} \subseteq
\gamma^iR$ for all positive integers $i$.   
\label{williams} 
\end{lem}

\begin{proof} This is \cite[Lemma 4.3]{Williams}.
\end{proof}

\begin{lem} Let $(R,m)$ be a normal local domain and $J$ an ideal of pure 
height one, which has order $n$ when regarded as an element of the divisor
class group $\Cl(R)$. Then for $0 < i < n$, we have $J^{(i)} J^{(n-i)}
\subseteq J^{(n)} m$. 
\label{reflexive}
\end{lem}

\begin{proof}  Let $J^{(n)} = \alpha R$. Clearly $J^{(i)} J^{(n-i)}\subseteq
\alpha R$, and it suffices to show that $J^{(i)} J^{(n-i)} \subsetneq \alpha
R$. If $J^{(i)} J^{(n-i)} = \alpha R$, then $J^{(i)}$ is an invertible
fractional ideal, and so must be a projective $R$--module. Since $R$ is local, 
$J^{(i)}$ is a free $R$--module, but this is a contradiction since $J^{(i)}$
cannot be principal for  $0 < i < n$. 
\end{proof}

\begin{disc} Let $(R,m)$ be a $\mathbb Q$--Gorenstein \CM normal local domain,
with canonical ideal $J$. Let $n$ denote the order of $J$ as an element of the
divisor class group $\Cl(R)$, and pick $\alpha \in R$ such that $J^{(n)} =
\alpha R$. Consider the subring $R[JT, J^{(2)}T^2, \dots]$ of $R[T]$, and
let 
$$  
S = R[JT, J^{(2)}T^2, \dots] /  (\alpha T^n -1).
$$
Note that $S$ has a natural 
${\mathbb Z} / n{\mathbb Z}$--grading where $[S]_0 = R$,
and for $0 < i < n$ we have  $[S]_i = J^{(i)}T^i$. We claim that the ideal 
$$
\mathfrak{m} = m + JT + J^{(2)}T^2 + \dots + J^{(n-1)}T^{n-1}
$$ 
is a maximal  ideal of $S$. Since each $J^{(i)}$ is an ideal of $R$, we
need only  verify that $J^{(i)}T^i \mathfrak{m} \subseteq \mathfrak{m}$ for $0
< i < n-1$, but this follows from Lemma \ref{reflexive}. Note furthermore that 
$\mathfrak{m}^n \subseteq mS$. \end{disc}

\section{The main result}

\begin{proof}[Proof of Theorem \ref{main}] 
The property of being a splinter localizes, as does the property of being
$\mathbb{Q}$--Gorenstein. Hence if the splinter ring $R$ is not
F--regular, we may
localize at a prime ideal $P \in \Spec R$ which is minimal with respect to the
property that $R_P$ is not F--regular. After a change of notation, we have a
splinter $(R,m)$ which has an isolated non F--regular point at the maximal
ideal $m$. This shows that $R$ has an $m$--primary test ideal. However since
$R$ is a splinter it must be F--pure, and so the test ideal is
precisely the maximal ideal $m$. Note that by 
\cite[Theorem 5.1]{Sminv} parameter ideals of $R$ are tightly
closed, and $R$ is indeed F--rational.

Let $\dim R = d$. Choose a system of parameters for $R$ as follows: first pick
a nonzero element $y_1 \in J$. Then, by Lemma \ref{williams}, pick $y_2 $ not
in any minimal prime of $y_1$ such that $y_2^i J^{(i)} \subseteq \gamma^i R$
for a fixed element $\gamma \in J$, for all positive integers $i$. Extend $y_1,
y_2$ to a full system of parameters $y_1, \dots, y_d$ for $R$. Since $y_1 \in
J$, there exists $u \in R$ such that $s=uy_1$ is a socle generator in $J /
(y_1, \dots, y_d)J$. Let $Y$ denote the product $y_1 \dots  y_d$.

Consider the family of ideals $\{I_c\}_{c \in \mathbb N}$ as in Lemma
\ref{irreducible}. If $R$ is not F--regular, there exists an irreducible ideal 
$I_c = (y_1^c, \dots, y_d^c)J :_R s$ which is not tightly closed, specifically 
$Y^{c-1} \in I_c^*$. Consequently $sY^{c-1} \in (y_1^c, \dots, y_d^c)J^*$.
In the ring $S$, this says that  $sY^{c-1} \in (y_1^c, \dots,
y_d^c)JS^*$ and so
$$
sTY^{c-1} \in (y_1^c, \dots, y_d^c)JTS^* \subseteq (y_1^c, \dots, y_d^c)S^*.
$$
We shall first imitate the proof of \cite[Lemma 5.2]{Sminv} to obtain an 
\lq\lq equational condition\rq\rq \ from this statement. To simplify notation,
let $z=sTY^{c-1}$ and $x_i = y_i^c$ for $1 \le i \le d$. We then have $z \in
(x_1, \dots, x_d)S^*$. Consider the maximal ideal $\mathfrak{m} = m + JT +
J^{(2)}T^2 + \dots + J^{(n-1)}T^{n-1}$ of $S$ and the highest local cohomology
module 
$$
H^{d}_\mathfrak{m}(S) = \varinjlim S/(x_1^i,\dots, x_d^i),
$$
where the maps in the direct limit system are induced by
multiplication by $x_1 \dotsm x_d$. 

Since the test ideal of $R$ is $m$, if $Q_0$ is a power of $p$ greater 
than $n$, we have 
$\mathfrak{m}^{Q_0} z^q \in (x_1^q, \dots, x_d^q)S$ for all $q = p^e$.

Let $\eta$ denote $[z + (x_1, \dots, x_d)S]$ viewed as an element of
$H^{d}_\mathfrak{m}(S)$, and $N$ be the $S$--submodule of 
$H^{d}_\mathfrak{m}(S)$ spanned by all $F^e(\eta)$ where  $e \in \mathbb{N}$.
Since $H^{d}_\mathfrak{m}(S)$ is an $S$--module with DCC, there exists $e_0$
such that the submodules generated by $F^{e_0}(N)$ and $F^{e'}(N)$ agree for
all $e' \ge e_0$. Hence there exists an equation of the form
$$ 
F^{e_0}(\eta) = a_1 F^{e_1}(\eta) + \dots + a_k F^{e_k}(\eta) 
$$ 
with $a_1, \dots, a_k \in S$ and $e_0 < e_1 \le e_2 \le \dots \le e_k$. If some
$a_i$ is not a unit, we may use suitably high Frobenius iterations on the
equation above, and the fact that for $Q_0 \ge n$ we have 
$\mathfrak{m}^{Q_0} F^e(\eta) = 0$ for all $e \in \mathbb N$, to replace
the above
equation by one in which the coefficients which occur are indeed units. 
Hence we have an equation $F^{e}(\eta) = a_1 F^{e_1}(\eta) + \dots + a_k
F^{e_k}(\eta)$ where $e < e_1 \le e_2 \le \dots \le e_k$ and  $a_1, \dots, a_k$
are units. Let $q = p^e$, $q_i = p^{e_i}$ for $1 \le i \le k$ and $X = x_1 \dots
x_d$.  Rewriting our equation we have
\begin{eqnarray*}
[z^q X^{q_k-q} + (x_1^{q_k}, \dots, x_d^{q_k})S] \ = \
     a_1 [z^{q_1} X^{q_k-q_1} + (x_1^{q_k}, \dots, x_d^{q_k})S]  
     \ + \ \cdots \\
       \cdots \ + \ a_k[z^{q_k} + (x_1^{q_k}, \dots, x_d^{q_k} )S],
\end{eqnarray*}
i.e., \ $[z^q X^{q_k-q} \ - \ a_1 z^{q_1} X^{q_k-q_1} \ - \ \cdots \ - 
             \ a_k z^{q_k} \ + \ (x_1^{q_k}, \dots, x_d^{q_k} )S] = 0$. 
Since the ring $S$ may not necessarily be Cohen--Macaulay, we cannot
assume that the maps in the direct limit system 
$\varinjlim S/(x_1^i,\dots, x_d^i)$ are injective. However for a
suitable positive integer $b$ we do obtain the equation
\begin{eqnarray*}
(zX^{b-1})^Q \ \in (x_1^{bQ},  \dots, \ x_d^{bQ}, \ zX^{bQ-1}, \
            z^p X^{bQ-p},  \dots, \ z^{Q/p} X^{bQ-Q/p} )S,
\end{eqnarray*}
where $Q=q_k$.
Going back to the earlier notation and setting $t=bc$, we have
$$
(sTY^{t-1})^Q \ \in (y_1^{tQ}, \dots, y_d^{tQ}, \ sT Y^{tQ-1}, \ (sT)^p Y^{tQ-p},
                  \dots, (sT)^{Q/p} Y^{tQ-Q/p} )S.
$$   
Note that $\frac{1}{T} = \alpha T^{n-1} \in S$, and multiplying the above by 
$\frac{1}{T^Q}$, we get
\begin{eqnarray*}
(sY^{t-1})^Q \ \in \ (y_1^{tQ}\frac{1}{T^Q}, \dots, \ y_d^{tQ}\frac{1}{T^Q}, \ 
               s Y^{tQ-1}\frac{1}{T^{Q-1}}, \  
               s^p  Y^{tQ-p}\frac{1}{T^{Q-p}}, \dots \\ \dots,  
               s^{Q/p} Y^{tQ-Q/p}\frac{1}{T^{Q-Q/p}} )S.
\end{eqnarray*}
Since $(sY^{t-1})^Q \in [S]_0 = R$, we may intersect the ideal above with $R$
to obtain
\begin{eqnarray*}
(sY^{t-1})^Q \in (y_1^{tQ}J^{(Q)}, \dots, \ y_d^{tQ}J^{(Q)}, \ 
               s Y^{tQ-1}J^{(Q-1)}, \  s^p  Y^{tQ-p}J^{(Q-p)}, \dots \\ \dots, \ 
               s^{Q/p} Y^{tQ-Q/p}J^{(Q-Q/p)} )R.
\end{eqnarray*}
Replacing $s = uy_1$ above, we get
\begin{eqnarray*}
(uy_1 Y^{t-1})^Q \ \in \ (y_1^{tQ}J^{(Q)}, \ \dots, \ y_d^{tQ}J^{(Q)}, \ 
               (uy_1) Y^{tQ-1}J^{(Q-1)}, \\
               (uy_1)^p  Y^{tQ-p}J^{(Q-p)}, \ \dots, \ 
               (uy_1)^{Q/p} Y^{tQ-Q/p}J^{(Q-Q/p)} )R.
\end{eqnarray*}
Let $Z=\frac{Y}{y_1} = y_2 \dots y_d$. We then have
\begin{eqnarray*}
(uZ^{t-1})^Q y_1^{tQ} \ \in \ (y_1^{tQ} J^{(Q)}, \ y_2^{tQ}, \ \dots, \ y_d^{tQ}, \ 
               uy_1^{tQ} Z^{tQ-1}J^{(Q-1)}, \\  
               u^p  y_1^{tQ}  Z^{tQ-p}J^{(Q-p)}, \ \dots, \ 
               u^{Q/p}  y_1^{tQ} Z^{tQ-Q/p}J^{(Q-Q/p)} )R.
\end{eqnarray*}
Using the fact that $y_1, \dots, y_d$ are a system of parameters for the \CM
ring $R$, we get
\begin{eqnarray*}
(uZ^{t-1})^Q \ \in \ (J^{(Q)}, \ y_2^{tQ}, \ \dots, \ y_d^{tQ}, \ 
               u Z^{tQ-1}J^{(Q-1)}, \  
               u^p  Z^{tQ-p}J^{(Q-p)}, \dots \\   
               \dots, \ u^{Q/p} Z^{tQ-Q/p}J^{(Q-Q/p)} )R.
\end{eqnarray*}
Consequently there exists $a \in J^{(Q)}$, $b_i \in R$ and $c_{p^e} \in
J^{(Q-Q/p^e)}$ such that
\begin{eqnarray*}
(uZ^{t-1})^Q \ = \  a \ + \ \sum_{i=2}^{d} b_i y_i^{tQ} \ + \ c_1 u Z^{tQ-1} \ + 
                \ c_p u^p  Z^{tQ-p} \ + \ \cdots \\
                \cdots \ + \ c_{Q/p} u^{Q/p} Z^{tQ-Q/p}.
\end{eqnarray*}

For $2 \le i \le d$, consider the following equations in the variables 
$V_2, \dots, V_d$:
\begin{eqnarray*}
V_i^Q  =  b_i + c_1 V_i (\frac{Z}{y_i})^{tQ-t} \ +  
            c_p V_i^p  (\frac{Z}{y_i})^{tQ-tp}  + \cdots 
            + c_{Q/p} V_i^{Q/p} (\frac{Z}{y_i})^{tQ-tQ/p}.
\end{eqnarray*}
Since these are monic equations defined over $R$, there exists a module finite
normal extension ring $R_1$, with solutions $v_i$ of these equations. Working 
in the ring $R_1$, let 
$$
w = uZ^{t-1}-\sum_{i=2}^{d} v_i y_i^t.
$$
Combining the earlier equations, we have
$$
w^Q = a + c_1 w Z^{tQ-t} + c_p w^p  Z^{tQ-tp} + \dots 
        + c_{Q/p} w^{Q/p} Z^{tQ-tQ/p}.
$$
Multiplying this equation by $y_2^Q$ and using the fact that 
$y_2^i J^{(i)} \subseteq \gamma^i R$ for all positive integers $i$, we get
$$
(wy_2)^Q = d_0 \gamma^Q + d_1 wy_2 \gamma^{Q-1} + d_p(wy_2)^p \gamma^{Q-p}
         + \dots + d_{Q/p} (wy_2)^{Q/p} \gamma^{Q-Q/p}.
$$
Since the ring $R_1$ is normal, 
The above equation gives an equation by which $wy_2 / \gamma$ is integral
over the ring $R_1$. Since $R_1$ is normal, we have $wy_2 \in \gamma R_1$.
Combining this with 
$w = uZ^{t-1}-\sum_{i=2}^d v_i y_i^t$, we have
$$
uZ^{t-1}y_2 = wy_2 + (\sum_{i=2}^d v_i y_i^t)y_2
  \ \in \ (J,\  y_2^{t+1}, \ y_2y_3^t, \ \dots, \ y_2y_d^t )R_1,
$$
and so 
$$
uZ^{t-1}y_2 \ \in \ (J, \ y_2^{t+1}, \ y_2 y_3^t,  \dots, \ y_2 y_d^t )^+  = 
(J, \ y_2^{t+1}, \ y_2 y_3^t,  \dots, \ y_2 y_d^t )R.
$$
Since $y_2$ is not in any minimal prime of $J$, we get $uZ^{t-1} \in 
(J,  y_2^t, y_3^t, \dots, y_d^t )R$.  Multiplying this by $y_1$, we get
$$ 
sZ^{t-1} \ \in \ (y_1 J, \ y_1 y_2^t, \ y_1 y_3^t, \ \dots, \ y_1 y_d^t )R 
         \subseteq (y_1, \ y_2^t, \ y_3^t, \ \dots, \ y_d^t)J,
$$
but this contradicts the fact that $s$ generates the socle in 
$J / (y_1, \dots, y_d)J$.
\end{proof}

\begin{cor} Let $(R,m)$ be an excellent integral domain of dimension
two over a field of characteristic $p > 0$. Then 
$R$ is a splinter if and only if it is F--regular.
\end{cor}

\begin{proof} The hypotheses imply that $R$ is F--rational, and so has a
torsion divisor class group by a result of Lipman, \cite{Li}. Hence $R$
must be $\mathbb Q$--Gorenstein.
\end{proof}

\begin{defn} Let $R=K[X_1, \dots, X_n]/I$ be a domain finitely generated over a
field $K$ of characteristic zero. We say $R$ is of {\it splinter type}\/ if
there exists a  finitely generated $\mathbb Z$--algebra $A \subseteq K$ and a
finitely generated free $A$--algebra $R_A = A[X_1, \dots, X_n]/I_A$ such that
$R \cong R_A \otimes_A K$, and for all maximal ideals $\mu$ in a Zariski dense
subset of $\Spec A$, the fiber rings $R_A \otimes_A A/\mu$ (which are rings
over fields of characteristic $p$) are splinter.
\end{defn}
 
Using the equivalence of F--regular type and log--terminal singularities for 
rings finitely generated over a field of characteristic zero (see
\cite{Hara, Smvanish, Walog}) we obtain the following corollary: 

\begin{cor} Let $R$ be a finitely generated $\mathbb Q$--Gorenstein domain over 
a field of characteristic zero. Then $R$ has log--terminal singularities if and
only if it is of splinter type.
\label{splinter-type}
\end{cor}

\section*{Acknowledgement}
 
It is a pleasure to thank Melvin Hochster for several valuable discussions on
tight closure theory.

\end{document}